\documentclass[11pt]{amsart}
\usepackage{amssymb, color}
\usepackage{amsmath}
\usepackage{amsthm}

\makeatletter
\def\LaTeX{\leavevmode L\raise.42ex
    \hbox{\kern-.3em\size{\sf@size}{0pt}\selectfont A}\kern-.15em\TeX}
\makeatother

\newcommand{\BibTeX}{{\rm B\kern-.05em{\sc
          i\kern-.025emb}\kern-.08em\TeX}}

\makeatletter
\def\@currentlabel{2.1}\label{e:dispaa}
\def\@currentlabel{2.21}\label{e:dispau}
\def\@currentlabel{2.22}\label{e:dispav}
\def\@currentlabel{2.23}\label{e:dispaw}
\def\@currentlabel{2.24}\label{e:dispax}
\def\theequation{\thesection.\@arabic\c@equation}
\makeatother

\renewcommand{\theequation}{\arabic{section}.\arabic{equation}}

\newtheorem{thm}{Theorem}[section]
\newtheorem{lem}[thm]{Lemma}
\newtheorem{cor}[thm]{Corollary}
\newtheorem{prop}[thm]{Proposition}
\theoremstyle{definition}

\newtheorem{rem}[thm]{Remark}

\newcommand{\R}{\mathbb{R}}
\newcommand{\N}{\mathbb{N}}
\newcommand{\e}{\epsilon}

\def \D{\Delta}

\def \l{\lambda}
\def \O{\Omega}
\def \p{\partial}

\begin{document}
\title{On stable solutions of biharmonic problem with polynomial growth}

\author{Hatem Hajlaoui}
\address{Institut de math\'ematiques appliqu\'ees et d'informatiques, Kairouan, Tunisia}
\email{hajlaouihatem@gmail.com}

\author{Abdelaziz Harrabi}
\address{Institut de math\'ematiques appliqu\'ees et d'informatiques, Kairouan, Tunisia}
\email{abdellaziz.harrabi@yahoo.fr}

\author{Dong Ye}
\address{LMAM, UMR 7122,
Universit\'e de Lorraine-Metz,
57045 Metz, France}
\email{dong.ye@univ-lorraine.fr}



\date{}
\def\baselinestretch{1}

\begin{abstract}
We prove the nonexistence of smooth stable solution to the biharmonic
problem $\Delta^2 u= u^p$, $u>0$ in $\R^N$ for $1 < p < \infty$ and
$N < 2(1 + x_0)$, where $x_0$ is the largest root of the following
equation:
$$x^4 - \frac{32p(p+1)}{(p-1)^2}x^2 + \frac{32p(p+1)(p+3)}{(p-1)^3}x -\frac{64p(p+1)^2}{(p-1)^4} = 0.$$
In particular, as $x_0 > 5$ when $p > 1$, we obtain the
nonexistence of smooth stable solution for any $N \leq 12$ and $p > 1$.
Moreover, we consider also the corresponding problem in the half space $\R^N_+$, or the elliptic problem $\Delta^2 u=
\l(u+1)^p$ on a bounded smooth domain $\O$ with the Navier boundary
conditions. We will prove the regularity
of the extremal solution in lower dimensions. Our results improve
the previous works in \cite{wy, wxy, co, cg}.
\end{abstract}

\maketitle \baselineskip 16pt
\section{Introduction}
\setcounter{equation}{0}
Consider the biharmonic equation
\begin{equation}
\label{1.1} \Delta^2 u = u^p, \quad u>0 \;\; \mbox{in $\R^N$}
\end{equation}
where $N \geq 5$ and $p > 1$. Let
\begin{equation}
\label{1.2} \Lambda(\phi) := \int_{\R^N} |\Delta \phi|^2 dx -p
\int_{\R^N} u^{p-1} \phi^2 dx, \quad \forall\; \phi \in H^2(\R^N).
\end{equation}
A solution $u$ is said stable if $\Lambda(\phi) \geq 0$ for any test function $\phi \in H^2(\R^N)$.

\medskip
In this note, we prove the following classification result.
\begin{thm}
\label{main} Let $N \geq 5$ and $p > 1$. The equation \eqref{1.1}
has no classical stable solution, if $N < 2 + 2x_0$ with $x_0$ the
largest root of the polynomial
\begin{align}
\label{poly}
H(x) = x^4 - \frac{32p(p+1)}{(p-1)^2}x^2 + \frac{32p(p+1)(p+3)}{(p-1)^3}x -\frac{64p(p+1)^2}{(p-1)^4}.
\end{align}
Moreover, we have $x_0 > 5$ for any $p > 1$. Consequently, if $N
\leq 12$, \eqref{1.1} has no classical stable solution for all $p >
1$.
\end{thm}

\medskip
For the corresponding second order problem:
\begin{equation}
\label{1.4} \Delta u+|u|^{p-1} u=0 \ \ \mbox{in} \
\R^N, \quad p > 1.
\end{equation}
Farina has obtained the optimal Liouville type result for all finite Morse index solutions.
 He proved in \cite{F} that a smooth finite Morse index solution to (\ref{1.4}) exists, if and only if $p \geq p_{JL}$ and $N \geq 11$, or $p = \frac{N+2}{N-2}$ and $N \geq 3$. Here $p_{JL}$ is the so-called Joseph-Lundgren exponent, see \cite{gnw}.

\medskip
The nonexistence of positive solutions to (\ref{1.1}) are showed if
$p < \frac{N+4}{N-4}$, and all entire solutions are classified if $p
= \frac{N+4}{N-4}$, see \cite{L, WX}. On the other hand, the
radially symmetric solutions to (\ref{1.1}) are studied in
\cite{FGK, gg, gw, k}. In particular, Karageorgis proved that the radial entire solution to \eqref{1.1} is stable if and only if $p \geq p_{JL_4}$ and $N \geq 13$. Here $p_{JL_4}$ stands for the corresponding Joseph-Lundgren
exponent to $\Delta^2$, see \cite{k}.

\medskip
The general fourth order case \eqref{1.1} is more delicate, since the integration by parts argument used by Farina cannot be adapted easily. The first nonexistence result for general stable solution was proved by Wei and Ye in \cite{wy}, they proposed to consider \eqref{1.1} as a system
\begin{align}
\label{sys} -\Delta u = v, \ -\Delta v = u^p \;\; \mbox{in} \ \R^N,
\end{align}
and introduced the idea to use different test functions with $u$ but
also $v$. Using estimates in \cite{S} they showed that for $N \leq
8$, \eqref{1.1} has no smooth stable solutions. For $N \geq 9$,
using blow-up argument, they proved that the classification holds
still for $p < \frac{N}{N-8} + \e_N$ with $\e_N > 0$, but without
any explicit value of $\e_N$. This result was improved by Wei, Xu
\& Yang in \cite{wxy} for $N \geq 20$ with a more explicit bound.

\medskip
Very recently, using the stability for system \eqref{sys} and some interesting iteration argument, Cowan proved that,
 see Theorem 2 in \cite{co}, there is no smooth stable solution to \eqref{1.1}, if $N < 2 + \frac{4(p+1)}{p-1}t_0$, where
\begin{align}
\label{t0}
t_0 = \sqrt{\frac{2p}{p+1}} + \sqrt{\frac{2p}{p+1} - \sqrt{\frac{2p}{p+1}}}, \quad \forall\; p > 1.
\end{align}
In particular, if $N \leq 10$, \eqref{1.1} has no stable solution
for any $p > 1$.

\medskip
However, the study for radial solutions in \cite{k} suggests the following conjecture:
$$\mbox{\sl A smooth stable solution to \eqref{1.1} exists if and only if $p \geq p_{JL_4}$ and $N \geq 13$.}$$
Consequently, the Liouville type result for stable solutions of \eqref{1.1} should hold true for $N\leq 12$
with any $p > 1$, that's what we prove here. More precisely, by Theorem 1 in \cite{k}, the radial entire solutions to \eqref{1.1} are unstable if and only if
\begin{align}
\label{JL4} \frac{N^2(N-4)^2}{16} < pQ_4\left(-\frac{4}{p-1}\right),
\quad \mbox{where } Q_4(m) = m(m-2)(m+N-2)(m+N-4).
\end{align}
The l.h.s.~ comes from the best constant of the Hardy-Rellich
inequality (see \cite{re}): Let $N \geq 5$,
\begin{align*}
\int_{\R^N} |\Delta\varphi|^2 dx \geq
\frac{N^2(N-4)^2}{16}\int_{\R^N}\frac{\varphi^2}{|x|^4} dx, \quad
\forall \; \varphi \in H^2(\R^N).
\end{align*}
The r.h.s~of \eqref{JL4} comes from the weak radial solution $w(x) =
|x|^{-\frac{4}{p-1}}$. When $p > \frac{N+4}{N-4}$, we can check that $w \in H^2_{loc}(\R^N)$ and
\begin{align*}
\Delta^2 w = Q_4\left(-\frac{4}{p-1}\right)w^p \ \ \mbox{in } {\mathcal D}'(\R^N).
\end{align*}
As $w^{p-1}(x) = |x|^{-4}$, using the Hardy-Rellich inequality, the condition \eqref{JL4} means just that $w$ is not a stable solution in $\R^N$, i.e.
$$\exists\; \varphi \in H^2(\R^N) \mbox{ such that } \Lambda_w(\varphi) := \int_{\R^N} |\Delta\varphi|^2 dx - p\int_{\R^N} Q_4\left(-\frac{4}{p-1}\right)w^{p-1}\varphi^2 dx < 0.$$
If we denote $N = 2 + 2x$, a direct calculation shows that
\eqref{JL4} is equivalent to $H_{JL_4}(x) < 0$, where
$$H_{JL_4}(x) := (x^2 - 1)^2 - \frac{32p(p+1)}{(p-1)^2}x^2 + \frac{32p(p+1)(p+3)}{(p-1)^3}x -\frac{64p(p+1)^2}{(p-1)^4}.$$
By \cite{gg}, \eqref{JL4} is equivalent to $N < 2 + 2x_1$ if $x_1$
denotes the largest root of $H_{JL_4}$. We can remark the nearness
between the polynomial $H$ in Theorem \ref{main} and $H_{JL_4}$, since $H(x) - H_{JL_4}(x) = 2x^2 - 1$.

\smallskip
Furthermore, Theorem \ref{main} improves the bound given in
\cite{co} for all $p > 1$. Indeed, there holds $x_0>
\frac{2(p+1)}{p-1}t_0$, see Lemmas \ref{L} and \ref{x0} below.

\medskip
Recall that to handle the equation \eqref{1.1}, we prove in general that $v = -\Delta u > 0$ in $\R^N$ using average function on the sphere, see \cite{WX}. Applying the blow up argument as in \cite{S, wy}, we can assume then $u$ and $v$ are uniformly bounded in $\R^N$. Therefore the following Souplet's estimate in \cite{S} holds true in $\R^N$, which was established for any {\it bounded} solution $u$ of \eqref{1.1}:
\begin{align}
\label{Sou}
v \geq \sqrt{\frac{2}{p+1}}u^\frac{p+1}{2}.
\end{align}

\medskip
Here we propose a new approach. Without assuming the boundedness of $u$ or showing immediately the positivity of $v$, we prove first some integral estimates for stable solutions of \eqref{1.1}, which will enable us the estimate \eqref{Sou}. This idea permits us to handle more general biharmonic equations: Let $N \geq 5$ and $p > 1$, consider
\begin{equation}
\label{1.1new} \Delta^2 u = u^p, \quad u>0 \;\; \mbox{in $\Sigma \subset \R^N$}, \quad u = \Delta u = 0\;\; \mbox{on $\p\Sigma$}.
\end{equation}
Let $E = H^2(\Sigma)\cap H_0^1(\Sigma)$ and
\begin{equation}
\label{1.2new} \Lambda_0(\phi) := \int_{\Sigma} |\Delta \phi|^2 dx -p
\int_{\Sigma} u^{p-1} \phi^2 dx, \quad \forall\; \phi \in E.
\end{equation}
A solution $u$ of \eqref{1.1new} is said stable if $\Lambda_0(\phi) \geq 0$ for any $\phi \in E$.
\begin{prop}
\label{newprop}
Let $u$ be a classical stable solution of \eqref{1.1new} with $\Sigma = \R^N$, or the half space $\Sigma = \R^N_+$, or the exterior domain $\Sigma = \R^N\setminus\overline{\Omega}, \R^N_+\setminus\overline{\Omega}$ where $\Omega$ is a bounded smooth domain of $\R^N$. Then the inequality \eqref{Sou} holds in $\Sigma$, consequently $v > 0$ in $\Sigma$.
\end{prop}
Using this, we obtain the Liouville type result for \eqref{1.1new} on half space situation, which improves the result in \cite{wy} for wider range of $N$, and without assuming the boundedness of $u$ or $v = -\Delta u$.
\begin{thm}
\label{newmain}
Let $x_0$ be defined in Theorem \ref{main}. If $N < 2 + 2x_0$, there exists no classical stable solution of \eqref{1.1new} if $\Sigma = \R^N_+$.
\end{thm}

\medskip
Our proof combines also many ideas coming from the previous works \cite{wy, cg, co}. Briefly, for \eqref{1.1}, we use different test functions to both equations of the system \eqref{sys} and we make use of the following inequality in \cite{cg} (see also \cite{co, dfs}): If $u$ is a stable solution of \eqref{1.1}, then
\begin{align}
\label{stable} \int_{\R^N} \sqrt{p}u^\frac{p-1}{2}\varphi^2 dx \leq
\int_{\R^N} |\nabla \varphi|^2 dx, \quad \forall\; \varphi \in
C_0^1(\R^N).
\end{align}
This will enable us two estimates. By suitable
combination, we prove that for any stable solution $u$ to
\eqref{1.1}, $\phi \in C^2_0(\R^N)$ and $s \geq 1$, there holds
\begin{align}
\label{est1} L(s) < 0 \ \ \Rightarrow \int_{\R^N}
u^pv^{s-1}\phi^2 dx \leq C\int_{\R^N}
v^s\left(|\Delta(\phi^2)| + |\nabla\phi|^2\right)
dx
\end{align}
Here $L$ is a polynomial of degree 4, see \eqref{defL} below, and the
constant $C$ depends only on $p$ and $s$.
Applying then the iteration argument of Cowan in \cite{co}, we show that $u \equiv 0$ if $N < 2 + 2x_0$, which is a contradiction, since $u$ is positive.

\medskip
Using similar ideas, we consider the elliptic equation on bounded domains:
$$
\left\{
\begin{array}{ll}
\Delta^2 u = \lambda(u + 1)^p & \mbox{ in a bounded smooth domain } \O \subset \R^N, \; \; N\geq 1\\
u = \D u = 0 & \mbox{ on } \p\O.
\end{array}
\right.
\leqno{(P_\lambda)}$$
It is well known (see \cite{bg, GGS}) that there exists a critical value $\lambda^* > 0$ depending on $p > 1$ and $\O$ such that
\begin{itemize}
\item If $\l \in (0, \l^*)$, $(P_\l)$ has a minimal and classical solution $u_\l$ which is stable;
\item If $\l = \l^*$, $u^* = \lim_{\l\to\l^*} u_\l$ is a weak solution to $(P_{\l^*})$, $u^*$ is called the {\it extremal solution}.
\item No solution of $(P_\l)$ exists whenever $\l > \l^*$.
\end{itemize}
In \cite{ceg, wy}, it was proved that if $1 < p < \left(\frac{N -
8}{N}\right)_+^{-1}$ or equally when $N < \frac{8p}{p-1}$, the extremal solution $u^*$ is smooth. Recently, Cowan \&
Ghoussoub improve the above result by showing that $u^*$ is smooth
if $N < 2 + \frac{4(p+1)}{p-1}t_0$ with $t_0$ in \eqref{t0},
so $u^*$ is smooth for any $p > 1$ when $N \leq 10$. Our result is
\begin{thm}
\label{main2} The extremal solution $u^*$ is smooth if $N< 2 + 2x_0$
with $x_0$ given by Theorem \ref{main}. In particular, $u^*$ is
smooth for any $p > 1$ if $N \leq 12$.
\end{thm}
Remark that our proof does not use the {\it a priori} estimate of $v = -\Delta u$ as in \cite{ceg, cg}.

\medskip
The paper is organized as follows. We prove some preliminary results and Proposition \ref{newprop}
in section 2. The proofs of Theorems \ref{main}, \ref{newmain} and \ref{main2} are given
respectively in section 3 and 4.

\section{Preliminaries}
\setcounter{equation}{0}
We show first how to obtain the estimate \eqref{Sou} for stable solutions of \eqref{1.1new}. Our idea is to use the stability condition \eqref{1.2new} to get some decay estimate for stable solutions of \eqref{1.1new}. In the following, we denote by $B_r$ the ball of center $0$ and radius $r > 0$.
\begin{lem}
\label{lemnew}
Let $u$ be a stable solution to \eqref{1.1new} and $v = -\Delta u$, there holds
\begin{align}
\label{2.0new}
\int_{\Sigma\cap B_R}\left(v^2 + u^{p+1}\right) dx \leq CR^{N - 4 - \frac{8}{p-1}}, \quad \forall\; R > 0.
\end{align}
\end{lem}

\medskip
\noindent{\bf Proof.} We proceed similarly as in Step 1 of the proof for Theorem 1.1 in \cite{wy}, but we do not assume here that $v > 0$ or $u$ is bounded in $\Sigma$. For any $\xi \in C^4(\Sigma)$ verifying $\xi = \Delta\xi = 0$ on $\p\Sigma$ and $\eta \in C_0^\infty(\R^N)$, we have
\begin{align}
\label{2.4new}
\begin{split}
\int_\Sigma (\Delta^2\xi) \xi\eta^2 dx = & \;\int_\Sigma \left[\D(\xi\eta)\right]^2 dx + \int_\Sigma \left[-4(\nabla\xi \cdot \nabla\eta)^2 + 2\xi\Delta\xi|\nabla \eta|^2\right]dx \\
& + \int_\Sigma \xi^2 \Big[2\nabla(\Delta\eta)\cdot \nabla\eta + (\D\eta)^2\Big] dx.
\end{split}
\end{align}
The proof is direct as for Lemma 2.3 in \cite{wy}, noticing just that in the integrations by parts, all boundary integration terms on $\p\Sigma$ vanish under the Navier conditions for $\xi$.

\smallskip
Take $\xi = u$, a solution of \eqref{1.1new} into \eqref{2.4new}, there holds
 \begin{align*}
 & \int_{\Sigma} [\Delta (u\eta)]^2 dx -\int_{\Sigma} u^{p+1}\eta^2 dx \\
  = & \; 4 \int_{\Sigma} (\nabla u \nabla \eta)^2 dx + 2\int_{\Sigma} uv |\nabla \eta|^2 dx -\int_{\Sigma} u^2 \Big[2\nabla(\Delta\eta)\cdot \nabla\eta + (\D\eta)^2\Big] dx
 \end{align*}
where $v = -\Delta u$. Using $\phi = u\eta$ in \eqref{1.2new}, we obtain easily
\begin{align}
\label{2.10new}
\begin{split}
& \int_{\Sigma} \Big[( \Delta (u\eta))^2 + u^{p+1}\eta^2\Big] dx \\
\leq & \; C_1 \int_\Sigma \Big[|\nabla u|^2|\nabla \eta|^2 + u^2|\nabla\left(\Delta \eta\right) \cdot \nabla \eta | + u^2(\D\eta)^2\Big] dx + C_2\int_{\Sigma} uv|\nabla \eta|^2dx.
\end{split}
\end{align}
Here and in the following, $C$ or $C_i$ denotes generic positive constants independent on $u$, which could be changed from one line to another. As $\Delta(u\eta)= 2\nabla u \cdot\nabla\eta + u\Delta\eta - v\eta$, from \eqref{2.10new}, we get
\begin{align}
\label{2.11new}
\begin{split}
  & \int_{\Sigma} \left[v^2\eta^2 + u^{p+1}\eta^2\right] dx \\
\leq & \; C_1\int_{\Sigma} \Big[|\nabla u|^2|\nabla \eta|^2 + u^2 |\nabla\left(\Delta \eta\right) \cdot \nabla \eta | + u^2(\D\eta)^2\Big] dx + C_2\int_{\Sigma} uv|\nabla \eta|^2dx.
\end{split}
\end{align}
On the other hand, as $u = 0$ on $\p\Sigma$,
\begin{align*}
2\int_{\Sigma}|\nabla u|^2|\nabla \eta|^2 dx & = \int_{\Sigma}\Delta(u^2)|\nabla \eta|^2 dx + 2\int_{\Sigma}uv|\nabla \eta|^2 dx\\
& = \int_{\Sigma} u^2 \Delta(|\nabla \eta|^2) dx + 2\int_{\Sigma}uv|\nabla \eta|^2 dx.
\end{align*}
Input this into \eqref{2.11new}, we can conclude that
\begin{align}
\label{2.12new}
\begin{split}
& \int_{\Sigma} \left[v^2\eta^2 + u^{p+1}\eta^2\right] dx\\
 \leq & \; C_1\int_{\Sigma} u^2 \Big[|\nabla\left(\Delta \eta\right) \cdot \nabla \eta | + (\D\eta)^2 + \left|\Delta(|\nabla \eta|^2)\right|\Big] dx + C_2\int_{\Sigma}uv|\nabla \eta|^2 dx.
\end{split}
\end{align}

\medskip
Take $\eta=\varphi^m$ with $m >2$ and $\varphi \in C_0^\infty(\R^N)$, $\varphi \geq 0$, it follows that
\begin{align*}
\int_{\Sigma} uv |\nabla \eta|^2 dx = &\; m^2 \int_{\Sigma} uv \varphi^{2(m-1)} |\nabla \varphi|^2 dx \\
\leq & \;\frac{1}{2C} \int_{\Sigma} (v \varphi^m)^2 dx + C\int_{\Sigma} u^2 \varphi^{2(m-2)} |\nabla \varphi|^4 dx.
 \end{align*}
Choose now $\varphi_0$ a cut-off function in $C_0^\infty(B_2)$ verifying $0 \leq \varphi_0 \leq 1$, $\varphi_0=1$ for $|x|<1$. Input the above inequality into (\ref{2.12new}) with $\varphi = \varphi_0(R^{-1}x)$ for $R > 0$, $\eta = \varphi^m$ and $m = \frac{2p+2}{p-1} > 2$, we arrive at
 \begin{equation}
 \label{2.13new}
 \begin{split}
 \int_{\Sigma} \left(v^2 + u^{p+1}\right)\varphi^{2m} dx & \leq \frac{C}{R^4}\int_{\Sigma} u^2\varphi^{2m-4} dx\\
  & \leq \frac{C}{R^4}\left(\int_{\Sigma} u^{p+1}\varphi^{(p+1)(m-2)} dx\right)^\frac{2}{p+1} R^{\frac{N(p-1)}{p+1}}\\
  & = \frac{C}{R^4}\left(\int_{\Sigma} u^{p+1}\varphi^{2m} dx\right)^\frac{2}{p+1} R^{\frac{N(p-1)}{p+1}}.
 \end{split}
 \end{equation}
 Hence
 \begin{align*}
 \int_{\Sigma} u^{p+1}\varphi^{2m} dx \leq CR^{N - \frac{4(p+1)}{p-1}}.
 \end{align*}
Combining with \eqref{2.13new}, as $\varphi^{2m} = 1$ for $x \in B_R := \{x \in \R^N, \;|x| \leq R\}$, \eqref{2.0new} is proved. \qed

\medskip\noindent
{\bf Proof of Proposition \ref{newprop}.} Let
\begin{align*}
\zeta = \beta u^\frac{p+1}{2} - v, \quad \mbox{where } \ \beta = \sqrt{\frac{2}{p+1}}.
\end{align*}
Then a direct computation shows that $\Delta \zeta \geq \beta^{-1}u^\frac{p-1}{2}\zeta$ in $\Sigma$. Consider $\zeta_+ := \max(\zeta, 0)$, there holds, for any $R > 0$
\begin{align}
\label{2.11}
\int_{\Sigma\cap B_R} |\nabla \zeta_+|^2 dx = -\int_{\Sigma\cap B_R} \zeta_+\Delta\zeta dx + \int_{\p(\Sigma\cap B_R)}\zeta_+ \frac{\p\zeta}{\p\nu} d\sigma \leq \int_{\Sigma\cap \p B_R}\zeta_+ \frac{\p\zeta}{\p\nu} d\sigma.
\end{align}
Here we used $\zeta_+\Delta\zeta \geq 0$ in $\Sigma$ and $\zeta = 0$ on $\p\Sigma$. Denote now $S^{N-1}$ the unit sphere in $\R^N$ and
\begin{align*}
e(r)= \int_{S^{N-1}\cap(r^{-1}\Sigma)} \zeta_+^2(r\sigma) d\sigma \quad \mbox{for $r > 0$}.
\end{align*}
Remark that $\exists\; R_0 > 0$ such that
\begin{align}
\label{2.12}
\int_{\Sigma\cap \p B_r}\zeta_+ \frac{\p\zeta}{\p\nu} d\sigma = \frac{r^{N-1}}{2}e'(r), \quad \forall\; r \geq R_0.
\end{align}
Moreover, for $R \geq R_0$, we deduce from \eqref{2.0new}
\begin{align*}
\int_{R_0}^R r^{N-1}e(r)dr \leq \int_{B_R\cap\Sigma} \zeta_+^2 dx \leq C\int_{B_R\cap\Sigma} \left(v^2 + u^{p+1}\right)dx \leq CR^{N - 4 - \frac{8}{p-1}} = o\left(R^N\right).
\end{align*}
This means that the function $e$ cannot be nondecreasing at infinity, so that there exists $R_j \to \infty$ satisfying $e'(R_j) \leq 0$. Combining with \eqref{2.11} and \eqref{2.12} with $R = R_j\to \infty$, there holds
\begin{align*}
\int_\Sigma |\nabla \zeta_+|^2 dx = 0.
\end{align*}
Using $\zeta = 0$ on $\p\Sigma$, we have $\zeta_+ \equiv 0$ in $\Sigma$, or equivalently \eqref{Sou} holds true in $\Sigma$. Clearly $v > 0$ in $\Sigma$ by \eqref{Sou}.\qed

\medskip
In the following, we show some properties of the polynomials $L$ and $H$, useful for our proofs. Let
\begin{align}
\label{defL}
L(s)=s^4-32\frac{p}{p+1}s^2+32\frac{p(p+3)}{(p+1)^2}s-64\frac{p}{(p+1)^2}, \quad s \in \R.
\end{align}
\begin{lem}
\label{L} $L(2t_0)<0$ and $L$ has a unique root $s_0$ in the interval
$(2t_0, \infty)$.
\end{lem}

\noindent
{\bf Proof.} Obviously
\begin{equation}
L(2t_0)=16t_0^4-128\frac{p}{p+1}t_0^2+64\frac{p(p+3)}{(p+1)^2}t_0-64\frac{p}{(p+1)^2}\nonumber
\end{equation}
By
$\frac{t_0^2}{2t_0-1}=\sqrt{\frac{2p}{p+1}}$ (see \cite{co}), there holds $t_0^4 = \frac{2p}{p+1}(2t_0 - 1)^2$. A direct computation
yields
\begin{align*}
    \frac {(p+1)^2L(2t_0)}{32p} & = (p+1)(2t_0 - 1)^2 - 4(p+1)t_0^2 + 2(p+3)t_0 - 2\\
    & =(p-1)(1-2t_0).
\end{align*}
As $t_0 > 1$ for any $p>1$, we have $L(2t_0)<0$.
Furthermore, $\forall\; p > 1$, $s\geq 2t_0$, we have
\begin{align*}
    (p+1)L''(s)=12(p+1)s^2-64p&\geq 48(p+1)t_0^2-64p\\
    & \geq 48(p+1)\frac{2p}{p+1} - 64p\\
    & =32p > 0
\end{align*}
 in $[2t_0, \infty)$, where we used $t_0^2 \geq \frac{2p}{p+1}$ by \eqref{t0}. Therefore $L$ is convex in $[2t_0, \infty)$, as $\lim_{s\to \infty} L(s) = \infty$ and $L(2t_0) < 0$, it's clear that $L$ admits a unique root in $(2t_0, \infty)$. \qed

\begin{rem}
Performing the change of variable $x=\frac{p+1}{p-1}s$, a direct calculation gives
\begin{align*}
H(x) = \left(\frac{p+1}{p-1}\right)^4L(s), \quad \mbox{hence } H(x)<0 \;\mbox{ if and only if }\; L(s)<0.
\end{align*}
Using the above Lemma, $x_0 = \frac{p+1}{p-1}s_0$ is the
largest root of the polynomial $H$, and $x_0$ is the unique root of $H$ for $x\geq \frac
{2(p+1)}{p-1}t_0$.
\end{rem}

\begin{lem}
\label{x0}
Let $x_0 = \frac{p+1}{p-1}s_0$ be the largest root of $H$. Then
$x_0>5$ for any $p>1$.
\end{lem}

\noindent
{\bf Proof}. As $x_0$ is the largest root of $H$, to have $x_0 > 5$, it is sufficient to show $H(5)<0$. Let $J(p) = (p-1)^4H(5)$, then $J(p)=-15p^4-1284p^3+4262p^2-3844p+625$. Therefore,
$$J'(p)=-60p^3-3852p^2+8524p-3844,\quad J''(p)=-180p^2-7704p + 8524.$$
We see that $J'' < 0$ in $[2, \infty)$. Consequently $J'(p) < 0$ and $J(p) < 0$ for $p \geq 2$. Hence $x_0> 5$ if $p \geq 2$. For $p \in (1, 2)$, there holds
$x_0> \frac{2(p+1)}{p-1}t_0 \geq 6t_0 > 5$ as $t_0 > 1$. \qed

\section{Proof of Theorems 1.1 and 1.3}
\setcounter{equation}{0}
We will prove only Theorem \ref{main}, since the proof of Theorem \ref{newmain} is completely similar, where we can change just $B_r$ by $B_r\cap\R^N_+$.

\medskip
The following result generalizes Lemma 4 in \cite{co}, which is a crucial argument for our proof. As above, the constant $C$ always denotes a
positive number which may change term by term, but does not depend on the solution $u$. For $k \in \N$, let $R_k:=2^kR$ with $R > 0$.
\begin{lem}
Assume that $u$ is a classical stable solution of (1.1). Then
for all $2\leq s<s_0$, there is $C <\infty$ such that
\begin{equation}
\label{estimate}
     \int_{B_{R_k}}u^{p}v^{s-1} dx \leq\frac{C}{R^2}\int_{B_{R_{k+1}}}v^{s} dx, \quad \forall\; R > 0.
\end{equation}
\end{lem}

\noindent
{\bf Proof}. Let $u$ be a classical stable solution of (1.1). Let $\phi \in C_0^2(\R^N)$ and $\varphi = u^{\frac{q+1}{2}}\phi$ with $q \geq 1$. Take $\varphi$ into the stability inequality \eqref{stable}, we obtain
\begin{equation}
\label{3.1}
\sqrt{p}\int_{\mathbb{R}^N}u^{\frac{p-1}{2}}u^{q+1}\phi^2\leq\int_{\mathbb{R}^N}u^{q+1}|\nabla\phi|^2
+\int_{\mathbb{R}^N}|\nabla u^{\frac{q+1}{2}}|^2\phi^2+ (q+1)\int_{\mathbb{R}^N}u^q\phi\nabla u\nabla\phi
\end{equation}
Integrating by parts, we get
\begin{align}
\label{3.2}
\begin{split}
\int_{\mathbb{R}^N}|\nabla
u^{\frac{q+1}{2}}|^2\phi^2 dx & = \frac{(q+1)^2}{4}\int_{\mathbb{R}^N}u^{q-1}|\nabla u|^2\phi^2 dx\\ & = \frac{(q+1)^2}{4q}\int_{\mathbb{R}^N}\phi^2\nabla(u^q)\nabla
u dx \\
&= \frac{(q+1)^2}{4q}\int_{\mathbb{R}^N}u^qv\phi^2 dx -\frac{q+1}{4q}\int_{\mathbb{R}^N}\nabla(u^{q+1})\nabla(\phi^2) dx \\
&= \frac{(q+1)^2}{4q}\int_{\mathbb{R}^N}u^qv\phi^2 dx + \frac{q+1}{4q}\int_{\mathbb{R}^N}u^{q+1}\Delta(\phi^2) dx
\end{split}
\end{align}
and
\begin{align}
\label{3.3}
  (q+1)\int_{\mathbb{R}^N}u^q\phi\nabla u\nabla\phi dx = \frac{1}{2}\int_{\mathbb{R}^N}\nabla(u^{q+1})\nabla(\phi^2) dx = -\frac{1}{2}\int_{\mathbb{R}^N}u^{q+1}\D(\phi^2) dx.
\end{align}
Combining \eqref{3.1}-\eqref{3.3}, we conclude that
\begin{align}
\label{3.4}
    a_1\int_{\mathbb{R}^N}u^{\frac{p-1}{2}}u^{q+1}\phi^2 dx \leq\int_{\mathbb{R}^N}u^qv\phi^2 dx +C\int_{\mathbb{R}^N}u^{q+1}\left(|\D(\phi^2)|+|\nabla
    \phi|^2\right) dx
\end{align}
where $a_1=\frac{4q\sqrt{p}}{(q+1)^2}$. Choose now $\phi(x)=
h(R_k^{-1}x)$ where $h \in C_0^\infty(B_2)$ such that $h \equiv 1$
in $B_1$, there holds then
\begin{align}
\label{3.5}
    \int_{\mathbb{R}^N}u^{\frac{p-1}{2}}u^{q+1}\phi^2 dx \leq\frac{1}{a_1}\int_{\mathbb{R}^N}u^qv\phi^2 dx + \frac{C}{R^2}\int_{B_{R_{k+1}}}u^{q+1} dx
\end{align}

Now, apply the stability inequality \eqref{stable} with
$\varphi = v^{\frac{r+1}{2}}\phi$, $r \geq 1$, there holds
\begin{align*}
\sqrt{p}\int_{\mathbb{R}^N}u^{\frac{p-1}{2}}v^{r+1}\phi^2\leq\int_{\mathbb{R}^N}v^{r+1}|\nabla\phi|^2
+\int_{\mathbb{R}^N}|\nabla v^{\frac{r+1}{2}}|^2\phi^2+ (r+1)\int_{\mathbb{R}^N}v^r\phi\nabla
v\nabla\phi
\end{align*}
By very similar computation as above (recalling that
$-\D v=u^p$), we have
\begin{align}
\label{3.7}
    \int_{\mathbb{R}^N}u^{\frac{p-1}{2}}v^{r+1}\phi^2 dx \leq \frac{1}{a_2}\int_{\mathbb{R}^N}u^pv^r\phi^2 dx +\frac{C}{R^2}\int_{B_{R_{k+1}}}v^{r+1} dx
\end{align}
where $a_2=\frac{4r\sqrt{p}}{(r+1)^2}$.

\medskip
Using \eqref{3.5} and \eqref{3.7}, there holds
\begin{align}
\label{3.8}
\begin{split}
I_1 + a_2^{r+1} I_2 & := \int_{\mathbb{R}^N}u^{\frac{p-1}{2}}u^{q+1}\phi^2 dx + a_2^{r+1}\int_{\mathbb{R}^N}u^{\frac{p-1}{2}}v^{r+1}\phi^2 dx\\ & \leq \frac{1}{a_1}\int_{\mathbb{R}^N} u^qv\phi^2 dx + a_2^{r}\int_{\mathbb{R}^N}
u^pv^r\phi^2 dx + \frac{C}{R^2}\int_{B_{R_{k+1}}}\left(u^{q+1} + v^{r+1}\right) dx.
\end{split}
\end{align}
Fix now
\begin{align}
\label{qr}
2q=(p+1)r+p-1, \quad \mbox{ or equivalently } \ q+1=\frac{(p+1)(r+1)}{2}.
\end{align}
By Young's inequality, we get
\begin{align*}
  \frac{1}{a_1}\int_{\mathbb{R}^N}u^qv\phi^2 dx & = \frac{1}{a_1}\int_{\mathbb{R}^N} {u^{\frac{p-1}{2}}u^{\frac{p+1}{2}r}v\phi^2} dx \\
  & = \frac{1}{a_1}\int_{\mathbb{R}^N} u^{\frac{p-1}{2}}u^{\frac{r}{r+1}(q+1)}v\phi^2 dx \\
  & \leq \frac{r}{r+1}
\int_{\mathbb{R}^N}u^{\frac{p-1}{2}}u^{q+1}\phi^2 dx +\dfrac
{1}{a_1^{r+1}(r+1)}\int_{\mathbb{R}^N}u^{\frac{p-1}{2}}v^{r+1}\phi^2 dx\\
& = \frac{r}{r+1}I_1 + \frac{1}{a_1^{r+1}(r+1)}I_2
\end{align*}
and similarly
\begin{align*}
  a_2^{r}\int_{\mathbb{R}^N} u^pv^r\phi^2 dx \leq \frac{1}{r+1} I_1 + \frac{a_2^{r+1}r}{r+1} I_2
\end{align*}
Combining the above two inequalities and \eqref{3.8}, we deduce then
\begin{align*}
a_2^{r+1}I_2\leq\left[\frac{a_2^{r+1}r}{r+1}+
\frac{1}{a_1^{r+1}(r+1)}\right]I_2 + \frac{C}{R^2}\int_{B_{R_{k+1}}}\left(u^{q+1} + v^{r+1}\right) dx,
\end{align*}
hence
\begin{align*}
    \frac{(a_1a_2)^{r+1}-1}{r+1} I_2\leq
    \frac{Ca_1^{r+1}}{R^2}\int_{B_{R_{k+1}}}\left(u^{q+1} + v^{r+1}\right) dx.
\end{align*}
Thus, if $a_1a_2 > 1$, by the choice of $\phi$,
\begin{align*}
    \int_{B_{R_{k}}}u^{\frac{p-1}{2}}v^{r+1} dx \leq I_2 \leq \frac{C}{R^2}\int_{B_{R_{k+1}}}\left(u^{q+1} + v^{r+1}\right) dx.
\end{align*}
From \eqref{Sou} and \eqref{qr}, we get $u^{q+1} \leq Cv^{r+1}$. Denote $s=r+1$, we can conclude that if $a_1a_2 > 1$,
\begin{align}
\label{3.10}
 \int_{B_{R_{k}}}u^pv^{s-1} dx \leq C_1\int_{B_{R_{k}}}u^\frac{p-1}{2}v^s dx \leq \frac{C_2}{R^2}\int_{B_{R_{k+1}}}\left(u^{q+1} + v^{r+1}\right) dx \leq \frac{C_3}{R^2}\int_{B_{R_{k+1}}}v^s dx.
\end{align}

On the other hand, a simple verification shows that
\begin{align*}
\mbox{$a_{1}a_2>1$ is equivalent to $L(s)<0$.}
\end{align*}
By Lemma \ref{L}, for $s \in [2t_0, s_0)$, there holds $L(s)<0$. So the inequality \eqref{3.10}, i.e.~\eqref{estimate}
holds true for any $2t_0 \leq s <s_0$. On the other hand, by Lemma 4 of \cite{co}, the estimate \eqref{estimate} is valid
for $2\leq s <2t_0$, hence for $2 \leq s < s_0$.\qed

\medskip
Now, we can follow exactly the iteration process in \cite{co} (see Proposition 1 or Corollary 2 there) to obtain
\begin{cor}
\label{co3.1}
Suppose $u$ is
a classical stable solution of (1.1).
 For all $2\leq \beta<\frac{N}{N-2}s_0$, there are $\ell\in \N$ and $C<\infty$ such that
 \begin{align*}
 \label{est2}
 \left(\int_{B_R}v^{\beta} dx\right)^\frac{1}{\beta}\leq
CR^{\frac{N}{2}(\frac{2}{\beta}-1)}\left(\int_{B_{R_{3\ell}}}v^2 dx\right)^{\frac{1}{2}}, \quad \forall\; R > 0.
\end{align*}
\end{cor}

Now we are in position to complete the proof of Theorem \ref{main}. Let $u$ be a
smooth stable solution to (1.1), combining Corollary \ref{co3.1}
and \eqref{2.0new}, for any $2\leq \beta < \frac{N}{N-2}s_0$, there
exists $C > 0$ such that
\begin{align*}
\left(\int_{B_R}v^{\beta} dx \right)^\frac{1}{\beta}\leq
CR^{\frac{N}{2}(\frac{2}{\beta}-1)+\frac{N}{2}-2-\frac{4}{p-1}}, \quad \forall\; R > 0.
\end{align*}
Note that
\begin{align*}
\frac{N}{2}\left(\frac{2}{\beta}-1\right)+\frac{N}{2}-2-\frac{4}{p-1} < 0\ \Leftrightarrow \ N<\frac{2(p+1)}{p-1}\beta.
\end{align*}
Considering the allowable range of $\beta$ given in Corollary \ref{co3.1}, if
$N<2+\frac{2(p+1)}{p-1}s_0$, after sending $R\rightarrow\infty$
we get then $\|v\|_{L^\beta(\R^N)} = 0$, which is impossible since $v$ is positive. To conclude, the equation \eqref{1.1} has no classical
stable solution if $N<2+2x_0$ where $x_0=\frac{p+1}{p-1}s_0$.

\medskip
Moreover, by Lemma \ref{x0}, $x_0>5$ for any $p>1$, which means that if $N\leq 12$, \eqref{1.1} has no classical stable solution
for all $p>1$. \qed

\section{Proof of Theorem 1.4}
\setcounter{equation}{0}
In this section, we consider the elliptic problem $(P_\l)$. Let $u_\l$ be the minimal solution of $(P_\l)$, it is well known that $u_\l$ is stable. To simplify the presentation, we erase the index $\l$. By \cite{cg, dfs}, there holds
\begin{equation}
\label{3.1new}
    \sqrt{\lambda p} \int_{\Omega}(u+1)^\frac{p-1}{2}\varphi^2 dx \leq\int_{\Omega}|\nabla\varphi|^2 dx, \quad \forall\; \varphi \in H_0^1(\Omega)
\end{equation}
Using $\varphi = u^\frac{q+1}{2}$ as test function in \eqref{3.1},
by similar computation as for \eqref{3.4} in section 3, we obtain
\begin{equation}
\label{4.1}
    a_1\sqrt{\lambda}\int_{\O}(u+1)^{\frac{p-1}{2}}u^{q+1} dx \leq\int_{\O}u^qv dx, \quad \mbox{where } \ a_1=\frac{4q\sqrt{p}}{(q+1)^2}.
\end{equation}
 Here we need not a cut-off function $\phi$, because all boundary terms appearing in the integrations by parts vanish under the Navier boundary conditions, hence the calculations are even easier. We can use the Young's inequality as for Theorem \ref{main}, but we show here a proof inspired by \cite{dggw}.

\medskip
Similarly as for \eqref{3.7}, using $\varphi = v^\frac{r+1}{2}$ in \eqref{3.1new}, we have
\begin{equation}
\label{4.2}
    a_2\sqrt{\lambda}\int_{\O}(u+1)^{\frac{p-1}{2}}v^{r+1} dx \leq\int_{\O}
    \lambda(u+1)^pv^r dx, \quad \mbox{where } \ a_2=\dfrac{4r\sqrt{p}}{(r+1)^2}.
\end{equation}
Take always $2q=(p+1)r+p-1$. Applying Holder's inequality, there hold
\begin{align}
\label{4.3}
\begin{split}
\int_{\O}u^qv dx & \leq
\left(\int_{\O}u^{\frac{p-1}{2}}v^{r+1} dx\right)^{\frac{1}{r+1}}\left(\int_{\O}u^{\frac{p-1}{2}+q+1} dx\right)^{\frac{r}{r+1}}\\
& \leq \left[\int_{\O}(u+1)^{\frac{p-1}{2}}v^{r+1} dx\right]^\frac{1}{r+1}\left(\int_{\O}u^{\frac{p-1}{2}+q+1} dx\right)^{\frac{r}{r+1}}
\end{split}
\end{align}and
\begin{equation}
\label{4.4}
\int_{\O}
(u+1)^pv^r dx \leq\left[\int_{\O}(u+1)^{\frac{p-1}{2}}v^{r+1} dx\right]^{\frac{r}{r+1}}\left[\int_{\O}(u+1)^{\frac{p-1}{2}+q+1} dx\right]^{\frac{1}{r+1}}.
\end{equation}
Multiplying \eqref{4.1} with \eqref{4.2}, using \eqref{4.3} and \eqref{4.4}, we get immediately
\begin{equation}
\label{4.5}
    \left[\int_{\O}(u+1)^{\frac{p-1}{2}}u^{q+1} dx\right]^{\frac{1}{r+1}}\leq
\frac{1}{a_1a_2}\left[\int_{\O}(u+1)^{\frac{p-1}{2}+q+1} dx\right]^{\frac{1}{r+1}}.
\end{equation}

On the other hand, for any $\varepsilon>0$ there exists $C_\varepsilon > 0$ such
that
\begin{align*}
(u+1)^{\frac{p-1}{2}+q+1}\leq
(1+\varepsilon)(u+1)^{\frac{p-1}{2}}u^{q+1}+C_\varepsilon \ \ \mbox{in } \ \R_+.
\end{align*}
If $a_1a_2 > 1$, there exists $\varepsilon_0 > 0$ satisfying $1 +
\varepsilon_0 < (a_1a_2)^{r+1}$. We deduce from \eqref{4.5} that
\begin{align*}
    \left[1- \frac{1+\varepsilon_0}{(a_1a_2)^{r+1}}\right] \int_{\O}(u+1)^{\frac{p-1}{2}}u^{q+1} dx \leq C.
\end{align*}
Therefore, when
$L(s)< 0$, i.e.~when $a_1a_2 > 1$, there is $C > 0$ such that
\begin{align*}
\int_\O u^{\frac{p-1}{2}+{q+1}} dx \leq \int_{\O}(u+1)^{\frac{p-1}{2}}u^{q+1} dx \leq C.
\end{align*}
As $u^* = \lim_{\l \to \l^*} u_\l$, we conclude, using Lemma \ref{L},
\begin{align}
\label{4.6}
u^* \in L^{\frac{p-1}{2}+{q+1}}(\O), \quad \mbox{for all $q$ satisifying } \frac{2(q+1)}{p+1} = r+1 = s < s_0.
\end{align}

Furthermore, by \cite{GGS}, we know that $u^* \in H^2(\O)$. As $u^* \geq 0$ verifies $\Delta^2 u^* = \l^*(u^*+1)^p \leq C(u^*)^{p-1} u^* + C$ with $u^* = \Delta u^* = 0$ on $\p\O$, by standard elliptic estimate, we know that $u^*$ is smooth if
\begin{align*}
  \frac{N}{4}  < \left(\frac{p-1}{2}+q+1\right)\frac{1}{p-1} = \frac{1}{2}\left(1+\dfrac{p+1}{p-1}s\right).
\end{align*} Therefore, $u^*$ is smooth if $N<2+2x_0$. By Lemma \ref{x0}, $u^*$ is smooth for any $p>1$
if $N\leq 12$.\qed

\medskip
\noindent {\bf Acknowledgments} {\it D.Y. is partially supported by the French
ANR project referenced ANR-08-BLAN-0335-01. This work was partially realized during a visit of A.H. at the University of Lorraine-Metz, he would like to thank {\sl Laboratoire de Math\'ematiques et Applications de Metz} for the kind hospitality.}

\end{document}